\documentclass[11pt,reqno]{amsart}
\usepackage[colorlinks=true,
pdfstartview=FitV, linkcolor=cyan, citecolor=magenta,
urlcolor=blue]{hyperref}
\usepackage{amsmath,amsfonts,latexsym,amssymb}
\usepackage{mathrsfs}
\usepackage[latin1]{inputenc}
\usepackage[T1]{fontenc}
\usepackage{ae,aecompl}
\usepackage{braket}
\usepackage{comment}
\usepackage{color}
\usepackage{graphicx}
\usepackage{subfig}
\newtheorem{theorem}{Theorem}[section]

\newtheorem{corollary}[theorem]{Corollary}

\long\def\@savemarbox#1#2{\global\setbox#1\vtop{\hsize\marginparwidth 
  \@parboxrestore\tiny\raggedright #2}}
\begin{document}
%%%%%%%%%%%%%%%%%%%%
\title[Hitchin representations]
{Hitchin representations of Fuchsian groups}
\author[Canary]{Richard Canary}
\address{University of Michigan, Ann Arbor, MI 41809}
\thanks{The author was partially supported by the grant DMS-1906441 from the National Science Foundation (NSF) and grant 674990
from the Simons Foundation}

\begin{abstract}
We survey the theory of Hitchin representations of closed surface groups into $\mathsf{PSL}(d,\mathbb R)$ with a focus on their
dynamical and geometric properties. We then describe recent extensions of this work to study Hitchin representations of co-finite area Fuchsian groups.
The motivation for this recent work is a conjecture about the geometry of the augmented Hitchin component.

\end{abstract}

\maketitle

\begin{verse}
{\em Dedicated to Dennis Sullivan on the occasion of his 80th birthday. Dennis was very kind to me when I was a feckless young
mathematician and he continues to be an inspiration now that I am a feckless old mathematician.}
\end{verse}

\section{Introduction}

Nigel Hitchin \cite{hitchin}  used the theory of Higgs bundles to enumerate the components of the ``character variety'' of (conjugacy classes of) representations
of a closed surface group into $\mathsf{PSL}(d,\mathbb R)$. He identified a component which is topologically a cell. When $d=2$, this component is the classical Teichm\"uller space
of a surface, so he called this component the Teichm\"uller component. This component is now known as the Hitchin component and representations in this
component are known as Hitchin representations.
In his paper, he makes a comment which served as one of the  primary motivations for
the new field of Higher (rank) Teichm\"uller theory.

\medskip

\begin{verse}
{\em  ``Unfortunately, the analytical point of view used for the proofs gives no indication of the
geometrical significance of the Teichm\"uller component.''}\\
--------Nigel Hitchin \cite{hitchin}
\end{verse}

\medskip

This challenge was taken up by Fran\c cois Labourie \cite{labourie-invent}, from a dynamical viewpoint, and Vladimir Fock and Alexander Goncharov \cite{fock-goncharov}, from
a more algebraic viewpoint. As one consequence, they were able to show that all Hitchin representations are discrete and faithful, and Labourie shows that they are quasi-isometric embeddings.
Subsequently, Andres Sambarino \cite{sambarino-quantitative} constructed Anosov flows which encoded spectral data associated to a Hitchin
representation. This allowed him to invoke results from the Thermodynamic Formalism to obtain counting and equidistribution results for Hitchin representations.
Martin Bridgeman, Dick Canary, Fran\c cois Labourie and Andres Sambarino \cite{BCLS,BCLS2}, built on Sambarino's earlier work and work of Curt McMullen \cite{mcmullen-wp}
in the classical setting to construct mapping class group invariant analytic Riemannian metrics on the Hitchin component
which generalize the Weil-Petersson metric on Teichm\"uller space.  These metrics are called pressure metrics.
Nicolas Tholozan \cite{deroin-tholozan} described an embedding of the Hitchin component into the Teichm\"uller space of foliated complex structures on
the unit tangent bundle of the surface, first studied by Dennis Sullivan \cite{sullivan-universal},  so that the pull-back of a ``Weil-Peterson metric'' on the Teichm\"uller space
is the simple root pressure metric.
(This is a highly selective history which reflects the focus of the author and hence of this paper.)

If one pinches a collection of disjoint curves in a hyperbolic surface, one naturally obtains a (possibly disconnected) cusped finite area hyperbolic surface. 
Bill Abikoff \cite{abikoff-aug} bordified
Teichm\"uller space by appending all such limiting surfaces and the result is known as augmented Teichm\"uller space. The ``strata at infinity'' are naturally
Teichm\"uller spaces of (possibly disconnected) finite area hyperbolic surfaces. This space itself is not well-behaved topologically, e.g. it is not locally compact,
but one may view it as the ``universal cover'' of the Deligne-Mumford compactification of moduli space. More explicitly, the quotient of augmented Teichm\"uller space
by the action of the mapping class group may be identified with the Deligne-Mumford compactification of Moduli space. Howard Masur \cite{masur-wp} showed
that the metric completion of  Teichm\"uller space with the Weil-Petersson metric may be identified with augmented Teichm\"uller space.

In recent years, the author has been fascinated by the goal of developing a theory of an  ``augmented Hitchin component'' and proving that it
arises as the metric completion of the Hitchin component with respect to a pressure metric.  The analogy is especially compelling when $d=3$, as
Hitchin representations are holonomy maps of convex (real) projective structures (see Choi-Goldman \cite{goldman-choi}) and points in the augmented
Hichin component conjecturally correspond to (possibly disconnected) finite area  convex projective surfaces. John Loftin and Tengren Zhang \cite{loftin-zhang} worked
out the topological aspects of the augmented Hitchin component when $d=3$ and obtain local  parametrizations of neighborhoods of the ``strata at infinity.''
In collaboration with Tengren Zhang and Andy Zimmer \cite{CZZ}, we  developed a theory of Hitchin representations of general geometrically finite Fuchsian groups.
In collaboration with Harry Bray, Nyima Kao and Giuseppe Martone \cite{BCKM}, we developed a dynamical framework to establish parallels of the counting
and equidistribution results of Sambarino. In subsequent work \cite{BCKM2}, we combine these results to construct pressure metrics on components of 
Hitchin representations of Fuchsian lattices, which arise as the strata at infinity for the augmented Hitchin component.

The first part of this paper will focus on the developments
of the classical Hitchin component, while the second portion will describe  a conjectural geometric picture of the augmented Hichin component and
describe the progress made towards this still elusive goal.

\medskip\noindent
{\bf Acknowledgements:} I would like to thank Fran\c cois Labourie, Giuseppe Martone, Andres Sambarino, Nicolas Tholozan and Tengren Zhang for helpful comments
on earlier versions of this manuscript. I also thank the referees for their careful reading of the original manuscript and suggestions which improved the exposition.

\section{The dynamical viewpoint of Labourie}

Fran\c cois Labourie \cite{labourie-invent} introduced techniques from Anosov dynamics to study the flow on the flat bundle associated to a Hitchin representation.
We recall some notation before describing his work. 

Throughout this paper $S$ will be a closed, orientable, connected surface of genus $g\ge 2$. Let $\mathcal \tau_d:\mathsf{PSL}(2,\mathbb R)\to\mathsf{PSL}(d,\mathbb R)$
be the irreducible representation, which is unique up to conjugation. A representation $\rho:\pi_1(S)\to \mathsf{PSL}(d,\mathbb R)$ is said to be {\em $d$-Fuchsian}
if it is (conjugate to) the result of post-composing a Fuchsian representation  $\rho_0:\pi_1(S)\to\mathsf{PSL}(2,\mathbb R)$ (i.e. a discrete, faithful representation) with the
irreducible representation $\tau_d$. A representation $\rho:\pi_1(S)\to\mathsf{PSL}(d,\mathbb R)$ is said to be a {\em Hitchin representation} if it can be continuously
deformed to a $d$-Fuchsian representation. One then defines the {\em Hitchin component} $\mathcal H_d(S)$ to be the space of $\mathsf{PGL}(d,\mathbb R)$-conjugacy classes of
Hitchin representations, i.e.
$$\mathcal H_d(S)\subset \mathrm{Hom}(\pi_1(S),\mathsf{PSL}(d,\mathbb R))/\mathsf{PGL}(d,\mathbb R).$$
Hitchin proved that $\mathcal H_d(S)$ is a cell.

\begin{theorem}{\rm (Hitchin \cite{hitchin})}
The Hitchin component $\mathcal H_d(S)$ is a real analytic manifold which is (real analytically) diffeomorphic to $\mathbb R^{(d^2-1)(2g-2)}$.
\end{theorem}

The {\em Fuchsian locus} in $\mathcal H_d(S)$  consists of (conjugacy classes of) $d$-Fuchsian representations and is an embedded copy of
the Teichm\"uller space $\mathcal T(S)$ of $S$.

We may identify $S=\mathbb H^2/\Gamma$ as a hyperbolic surface, where $\Gamma\subset\mathsf{PSL}(2,\mathbb R)$,  so $\Gamma$ is identified with $\pi_1(S)$.
The unit tangent bundle $T^1S$ of $S$ is then the quotient $T^1\mathbb H^2/\Gamma$ of the unit tangent bundle of $\mathbb H^2$. 
Hitchin observes that every Hitchin representation $\rho:\pi_1(S)\to\mathsf{PSL}(d,\mathbb R)$ lifts to a representation $\hat \rho:\pi_1(S)\to\mathsf{SL}(d,\mathbb R)$
(see Culler \cite{culler-lifting} for general criteria guaranteeing lfiting which apply in this case).
The {\em flat bundle}
associated to a representation $\rho:\pi_1(S)\to\mathsf{SL}(d,\mathbb R)$ is formed as
$$E_\rho=(T^1\mathbb H^2\times \mathbb R^d)/\Gamma$$
where the action on the first factor is the standard action of $\Gamma$ on $T^1\mathbb H^2$ and the action on $\mathbb  R^d$ is given by $\hat\rho(\Gamma)$.

The geodesic flow $\{\phi_t:T^1S\to T^1S\}_{t\in\mathbb R}$ lifts to the geodesic flow $\{\tilde\phi_t:T^1\mathbb H^2\to T^1\mathbb H^2\}_{t\in\mathbb R}$ and then extends to a flow on
$\{\tilde\psi_t\}_{t\in\mathbb R}$ on $T^1\mathbb H^2\times \mathbb R^d$ which acts trivially on the second factor, i.e. $\tilde\psi_t(\vec v,\vec w)=(\tilde\phi_t(\vec v),\vec w)$.
The flow $\{\tilde\psi_t\}$ then descends to a flow $\{\psi_t\}$ on $E_\rho$ which ``extends'' $\{\phi_t\}$.  
(If you intended to be confusing, you would say that $\{\psi_t\}$ is the flow parallel to the flat connection.)
%\marginpar{\tiny describe Hopf param?}

One way of stating Labourie's fundamental
dynamical result is that the flat bundle admits a splitting into line bundles with certain contraction properties.
We recall that a flow $\{\psi_t\}$ on a vector bundle $V$ over a compact base $B$ is {\em contracting} if given some (any) continuous family $\{||\cdot||_b\}_{b\in B}$ of norms
on the fibers, there exists $C,c>0$ so that 
$$||\psi_t(\vec v)||_{\phi_t(b)}\le Ce^{-ct}||\vec v||_b$$
for all $b\in B$, $t>0$ and $\vec v\in V_b$.

\begin{theorem}{\rm (Labourie \cite{labourie-invent})}
If $\rho\in\mathcal H_d(S)$, then $E_\rho$ admits a flow-invariant splitting 
$$E_\rho=L_1\oplus L_2\oplus\cdots\oplus L_d$$
into line bundles so that the flow induced by $\psi_t$ is contracting on $L_i\otimes L_j^*$ if $i>j$.
\end{theorem}

Notice that this splitting lifts to a flow-invariant, $\Gamma$-equivariant splitting 
$$T^1\mathbb H^2\times\mathbb R^d=\tilde{L_1} \oplus \tilde L_2\oplus\cdots\tilde L_d.$$
If $\gamma\in\Gamma$, let $\vec v_i$ be a non-trivial vector in $\tilde L_i$  lying over a tangent vector to the axis of $\gamma$. (Notice that, by flow invariance, the line $\tilde L_i$
is independent of which point on the axis you pick.)
Then, since the splitting is flow invariant and $\Gamma$-equivariant, $\vec v_i$ must be an eigenvector for $\hat\rho(\gamma)$. Let $\ell_i(\rho(\gamma))$ denote
the associated eigenvalue of  $\hat\rho(\gamma)$, i.e. $\hat\rho(\gamma)(\vec v_i)=\ell_i(\hat\rho(\gamma))\vec v_i$. The fact that $L_i\otimes L_j^*$ is contracting
if $i>j$ implies that $|\ell_i(\hat\rho(\gamma))|>|\ell_j(\hat\rho(\gamma))|$ if $i>j$. In particular,  $\rho(\gamma)$ is loxodromic and if 
$\lambda_i(\rho(\gamma))=|\ell_i(\hat\rho(\gamma))|$, then
$$\lambda_1(\rho(\gamma))>\lambda_2(\rho(\gamma))>\cdots >\lambda_d(\rho(\gamma)).$$
Moreover, since the flow is contracting we see that
there exists $C,c>0$ so that if $\ell(\gamma)$ denotes the translation length of $\gamma$ on $\mathbb H^2$, then
$$\frac{\lambda_i(\rho(\gamma))}{\lambda_{i+1}(\rho(\gamma))}\ge Ce^{c\ell(\gamma)}$$ if
$\gamma\in\Gamma$ and $1\le i\le d-1$.

Labourie's splitting also gives rise to a $\rho$-equivariant, H\"older continuous limit map $\xi_\rho:\partial\mathbb H^2\to \mathcal F_d$ where $\mathcal F_d$ is the space of 
$d$-dimensional flags.
%(In fact, in this setting the existence of such a map is equivalent to the splitting.)
If $x\ne y\in\mathbb H^2$, consider a flow line (a.k.a. geodesic) joining $y$ to $x$ and choose non-trivial vectors $\vec v_i\in \tilde L_i$ lying over a vector tangent to the
flow line. One then defines
$$\xi_\rho(x)=\left(\langle \vec v_1\rangle, \langle \vec v_1, \vec v_2 \rangle,\ldots,\langle \vec v_1,\vec v_2,\ldots \vec v_{d-1}\rangle\right)\qquad\text{and}\qquad
\xi_\rho(y)=\left(\langle \vec v_d\rangle, \langle \vec v_d, \vec v_{d-1} \rangle,\ldots,\langle \vec v_d,\vec v_{d-1},\ldots \vec v_{2}\rangle\right).$$
In particular, $\xi_\rho(x)$ is transverse to $\xi_\rho(y)$.
The contraction properties of the flow imply that $\xi_\rho(x)$ does not depend on $y$ and that $\xi_\rho$ is H\"older.

We summarize these properties below.

\begin{theorem}{\rm (Labourie \cite{labourie-invent})}
\label{Hitchin properties}
If $\rho\in\mathrm{Hom}(\pi_1(S),\mathsf{PSL}(d,\mathbb R))$ is a Hitchin representation, then
\begin{enumerate}
\item
There exists a H\"older continuous $\rho$-equivariant map $\xi_\rho:\partial\mathbb H^2\to \mathcal F_d$ so that if $x\ne y$,
then $\xi_\rho(x)$ is transverse to $\xi_\rho(y)$.
\item
if $\gamma\in\pi_1(S)$, then $\rho(\gamma)$ is loxodromic and $\xi_\rho(\gamma^+)$ is the attracting flag of $\rho(\gamma)$, where $\gamma^+$ is
the attracting fixed point of $\gamma$,
\item
There exists $C,c>0$ so that 
$$\frac{\lambda_i(\rho(\gamma))}{\lambda_{i+1}(\rho(\gamma))}\ge Ce^{c\ell(\gamma)}$$ if
$\gamma\in\Gamma$ and $1\le i\le d-1$.
\end{enumerate}
\end{theorem}

The properties above do not characterize Hitchin representations. For example, if one takes the direct product of a Fuchsian representation and
the trivial 1-dimensional representation, it will satisfy all the properties above. Such representations are known as Barbot representations, since
they were first studied by Thierry Barbot \cite{barbot}. However, Oliver Guichard \cite{guichard-char} extended Labourie's work to provide the following
characterization.

\begin{theorem}{\rm (Guichard \cite{guichard-char})}
If $\rho\in\mathrm{Hom}(\pi_1(S),\mathsf{PSL}(d,\mathbb R))$, then $\rho$ is a Hitchin representation if and only if there exists
a continuous $\rho$-equivariant map $\xi:\partial\mathbb H^2\to \mathbb P(\mathbb R^d)$ so that if $\{x_1,\ldots,x_d\}$ are distinct points in
$\partial\mathbb H^2$, then 
$$\xi(x_1)\oplus\cdots\oplus\xi(x_d)=\mathbb R^d.$$
\end{theorem}

Notice that is is clear that Barbot representations do not satisfy this characterization, since the image of any such map cannot span $\mathbb R^3$
in this case. (In fact, small deformations of Barbot representations also admit limit maps which fail to satisfy Guichard's criterion.)

Labourie's work also allows him to extend Fricke's Theorem to the setting of Hitchin components. We recall that the mapping class group $\mathrm{Mod}(S)$
of a closed orientable surface is the group of (isotopy classes of) orientation-preserving self-homeomorphisms of $S$.

\begin{theorem}{\rm (Labourie \cite{labourie-energy})}
The mapping class group $\mathrm{Mod}(S)$ acts properly discontinously on $\mathcal H_d(S)$.
\end{theorem}

\medskip\noindent
{\bf Remarks:} (1) Hitchin's work \cite{hitchin} establishes analogous results for representations into all split real Lie groups, but we will only
discuss the case of $\mathsf{PSL}(d,\mathbb R)$, which is the most studied case. Hitchin's work uses the theory of Higgs bundles. Due to
the author's woeful ignorance, we will not discuss any of the subsequent work on Hitchin representations from this more analytic viewpoint.

(2) In Labourie's seminal paper, he more generally defines Anosov representations of a hyperbolic group into any semi-simple Lie group.
This definition was explored more fully by Guichard-Wienhard \cite{guichard-wienhard} and later by Gu\'eritaud-Guichard-Kassel-Wienhard \cite{GGKW},
Kapovich-Leeb-Porti \cite{KLP}, Bochi-Potrie-Sambarino \cite{BPS}, Kassel-Potrie, Tsouvalas \cite{kostas} and others.
In particular, analogues of all the basic properties discussed above exist in this setting. For generalizations of Fricke's Theorem to this setting
see Guichard-Wienhard \cite[Cor. 5.4]{guichard-wienhard} or Canary \cite[Thm. 6.4]{canary-dynamics}.

The theory of Anosov representations has emerged as a central language in Higher Teichm\"uller theory. For those interested in me blathering on endlessly
about Anosov representations, lecture notes \cite{canary-informal} are available on my webpage.

\section{The positive viewpoint of Fock-Goncharov}

Vladimir Fock and Alexander Goncharov \cite{fock-goncharov} characterize Hitchin representations as those representations which admit a positive limit map.

We begin by recalling the definition of a positive map of a subset of the circle into $\mathcal F_d$. This definition relies on work of Lusztig \cite{lusztig}
and others on positivity in semi-simple Lie groups. Given an ordered basis $\mathcal B$ for $\mathbb R^d$, we say that a unipotent element $A\in\mathsf{SL}(d,\mathbb R)$ 
is  {\em totally positive} with respect to $\mathcal B$, if its matrix with respect to $\mathcal B$ is upper triangular and all its minors (which are not
forced to be 0 by the fact that the matrix is upper triangular) are strictly positive. The set  $\mathcal U_{>0}(\mathcal B)$ of  unipotent, totally
positive, upper triangular matrices with respect to $\mathcal B$ is a semi-group. 
An ordered  $k$-tuple $(F_1,F_2,\ldots, F_k)$ of distinct flags in $\mathcal F_d$ is positive with respect to an ordered basis $\mathcal B=(b_1,\ldots,b_d)$
for $\mathbb R^d$ if $b_i\in F_1^{(i)}\cap F_k^{(d-i+1)}$ for all $i$, and there exists $u_2,\ldots,u_{k-1}\in \mathcal U_{>0}(\mathcal B)$ so that
$F_i=u_{k-1}\cdots u_{i}F_k$ for all $i=2,\ldots,k-1$. Here $F^{(k)}$ denotes the $k$-dimensional component of a flag $F$.
If $X$ is a subset of $S^1$ then a map $\xi:X\to\mathcal F_d$ is {\em positive} if whenever $(x_1,\ldots, x_n)$ is a cyclically ordered subset of
distinct points in $X$, then $(\xi(x_1),\ldots,\xi(x_n))$ is positive with respect to some ordered basis. (In fact, it suffices to only consider $4$-tuples of points.)

If $\mathcal B=\{e_1,e_2\}$ is the standard basis for $\mathbb R^2$, then 
$$\mathcal U_{>0}(\mathcal B)=\left\{\begin{bmatrix}1 & a\\ 0 & 1\end{bmatrix}\ :\ a>0\right\}$$
and one may check that a $n$-tuple of distinct points in $\mathcal F_2=\mathbb P(\mathbb R^2)\cong S^1$ is positive with respect to $\mathcal B$
if and only if $x_1=[e_1]$, $x_n=[e_2]$ and the other points proceed monotonically in the counter-clockwise direction. More generally,
one may check that a map $\xi:X\to\mathcal F_2$ is positive if and only if it is monotonic. 

Vladimir Fock and Alexander Goncharov obtain the following characterization of Hitchin representations.

\begin{theorem}{\rm (Fock-Goncharov \cite{fock-goncharov})}
A representation $\rho:\pi_1(S)\to\mathsf{PSL}(d,\mathbb R)$ is a Hitchin representation if and only if there exists
a positive $\rho$-equivariant map $\xi:\partial\mathbb H^2\to \mathbb P(\mathbb R^d)$.
\end{theorem}

They develop further structure which allows them to explicitly parametrize Hitchin components using natural algebraic data. This
viewpoint was further developed by Francis Bonahon and Guillaume Dreyer \cite{bonahon-dreyer,bonahon-dreyer2}.

\medskip\noindent
{\bf Remark:} Olivier Guichard and Anna Wienhard \cite{gw-positivity} developed a more general notion of a $\Theta$-positivity for
real semi-simple Lie groups and characterize exactly which Lie groups admit such structures. Guichard, Labourie and Wienhard
\cite{GLW} recently proved that $\Theta$-positive representations are Anosov and that 
the space  of $\Theta$-positive representations into a given Lie group contains entire components of the space of reductive
representations.
This theory encompasses both
Hitchin representations into split real Lie groups and maximal representations into Hermitian Lie groups of tube type, as well
as certain representations into $\mathsf{SO}(p,q)$  (when $p\ne q$) and four other exceptional Lie groups. In the case of $\mathsf{SO}(p,q)$,
Jonas Beyrer and Beatrice Pozzetti \cite{beyrer-pozzetti} further show that the set of $\Theta$-positive representations is exactly a collection of components
of the representation variety.
The fundamental
conjecture is that $\Theta$-positive representations account for all components of  representation varieties of surface groups
into simple Lie groups which consist entirely of discrete faithful representations. (For a Higgs bundle-theoretic perspective see Bradlow-Collier-Garc\'ia-Prada-Gothen-Oliveira \cite{BCGO}.)

\section{Sambarino's geodesic flows}

Andres Sambarino \cite{sambarino-quantitative} defined a family of Anosov flows associated to a Hitchin representation which
record the spectral data of the representation. We first recall some linear algebra so that we can state his results.

Let 
$$\mathfrak{a}=\{\vec a\in\mathbb R^d\ :\ a_1+\ldots+a_d=0\}$$
be the standard Cartan algebra for $\mathsf{PSL}(d,\mathbb R)$.
The space $\mathfrak{a}^*$ of linear functionals on $\mathfrak{a}$ is generated by
the simple roots $\{\alpha_1,\ldots,\alpha_{d-1}\}$ where 
$$\alpha_i(\vec a)=a_i-a_{i+1}.$$
The standard positive Weyl Chamber $\mathfrak{a}^+$ is then the set where all the $\alpha_i$ are non-negative, i.e.
$$\mathfrak{a}^+=\{\vec a\in\mathfrak{a}\ :\ a_1\ge\cdots\ge a_d\}=\{\vec a\in\mathfrak{a}\ :\ \alpha_i(\vec a)\ge 0\ \forall\ i\}.$$
We will also be interested in the linear functionals given by the fundamental weights $\omega_k$ and the Hilbert length $\omega_H$ where
$$\omega_k(\vec a)=a_1+\ldots+a_k\qquad\text{and}\qquad\omega_H(\vec a)=a_1-a_d=\omega_1(\vec a)+\omega_{d-1}(\vec a).$$

The {\em Jordan projection} $\nu:\mathsf{PSL}(d,\mathbb R)\to\mathfrak{a}^+$ records the spectral data associated to an element
of $\mathsf{PSL}(d,\mathbb R)$.
If $A\in\mathsf{PSL}(d,\mathbb R)$ has
generalized eigenvalues with moduli 
$$\lambda_1(A)\ge\cdots\ge\lambda_d(A)$$
then
$$\nu(A)=\big(\log\lambda_1(A),\ldots,\log\lambda_d(A)\big).$$

If $\rho\in\mathcal H_d(S)$, the {\em Benoist limit cone} $\mathcal B(\rho)$ encodes the spectral data of $\rho(\Gamma)$.
Explicitly,
$$\mathcal B(\rho)=\overline{\bigcup_{\gamma\in\Gamma}\mathbb R_+\nu(\rho(\gamma))}\subset\mathfrak{a}^+.$$
We will be interested in the collection  $\mathcal B^+(\rho)$ of linear functionals which are positive on $\mathcal B(\rho)\setminus \{\vec 0\}$, given
by
$$\mathcal B^+(\rho)=\left\{\phi\in\mathfrak{a}^*\ |\ \phi(\vec a)>0\ \forall\ \vec a\in\mathcal B(\rho)\setminus\{\vec 0\} \right\}.$$
Theorem \ref{Hitchin properties} implies that  $\mathcal B(\rho)\setminus\{\vec 0\}$ is contained in the interior of $\mathfrak{a}^+$, so
$$\Delta=\left\{\phi\in\mathfrak{a}^*\ |\ \phi=\sum c_i\alpha_i,\ c_i\ge 0\ \forall\ i,\ \text{ and }\sum c_i>0\right\}\subset\mathcal B(\rho)^+.$$
%\marginpar{\tiny notate $\Delta_d$?}
Notice that, in particular, each fundamental weight $\omega_k\in\Delta$ and $\omega_H\in\Delta$.

If $\phi\in \mathcal B^+(\rho)$, then one obtains a natural associated length function, given by
$$\ell^\phi(\rho(\gamma))=\phi(\nu(\gamma)).$$
If $d=2$ and $\phi=\alpha_1=\omega_H=2\omega_1$, then $\ell^\phi(\rho(\gamma))$ is just the usual
translation length $\ell(\rho(\gamma))$ of $\rho(\gamma)$. One may then consider an associated {\em $\phi$-topological entropy} 
$h^\phi(\rho)$ which records the exponential growth
rate of the number of (conjugacy classes of) elements of  $\phi$-length at most $T$. Concretely, let
$$R_T^\phi(\rho)=\{[\gamma]\in[\Gamma]\ :\ \ell^\phi(\rho(\gamma))\le T\} \qquad\text{and}\qquad  h^\phi(\rho)=\lim_{T\to\infty} \frac{\log\#R_T^\phi(\rho)}{T}$$
where $[\Gamma]$ is the set of conjugacy classes of elements of $\Gamma$.

We can now summarize some of Sambarino's work. We recall that  a flow space ${\bf U}_1$ is said to be H\"older orbit equivalent to a flow space ${\bf U}_2$
if there is a H\"older homeomorphism $f:{\bf U}_1\to {\bf U}_2$ which takes flow lines to flow lines (but does not necessarily preserve the time parameter).

\begin{theorem}{\rm (Sambarino \cite{sambarino-quantitative,sambarino-hyperconvex})}
If $\rho\in\mathcal H_d(S)$ and $\phi\in\mathcal B^+(\rho)$, then there exists an Anosov flow ${\bf U}^\phi(\rho)$ which is H\"older orbit equivalent to the
geodesic flow on $T^1(S)$ so that the period of the orbit of ${\bf U}^\phi_\rho$ associated to $[\gamma]\in[\Gamma]$ is given by $\ell^\phi(\rho(\gamma))$.
Moreover, the topological entropy of ${\bf U}^\phi(\rho)$ is exactly $h^\phi(\rho)$ and
$$\#R_T(\phi)\sim \frac{e^{h^\phi(\rho)T}}{h^\phi(\rho)T}, \text{ i.e. } \lim_{T\to\infty} \frac{\#(R_T(\phi))h^\phi(\rho)T}{e^{h^\phi(\rho)T}}=1.$$
\end{theorem}

If $\alpha_{i,j}\in\mathfrak{a}^*$ is given by $\alpha_{i,j}(\vec a)=a_i-a_j$ (so $\alpha_i=\alpha_{i,i+1}$) and $i>j$, then one may obtain ${\bf U}^{\alpha_{i,j}}(\rho)$ from
the contracting line bundle $L_i\otimes L_j^*$.  The flows ${\bf U}^{\alpha_i}(\rho)$ are known as the {\em simple root flows}. 
Similarly, if $\omega_1\in\mathfrak{a}^*$ is the first fundamental weight, i.e. $\omega_1(\vec a)=a_1$, 
then one may observe that $L_1$ is contracting and obtain
${\bf U}^{\omega_1}(\rho)$, which we call the {\em spectral radius flow}, in the same manner.

We now explain how to obtain a flow space ${\bf U}_L$, H\"older orbit equivalent to $T^1S$, from a contracting (H\"older) line bundle $L$ over $T^1S$.
One first lifts $L$ to the contracting line bundle $\tilde L$ over $T^1\mathbb H^2$ and considers
the associate principal $\mathbb R$-bundle $\widehat L$ over $T^1\mathbb H^2$ so that the fiber over $\vec v\in T^1\mathbb H^2$ is given by 
$\big(\tilde L|_{\vec v}\setminus \{\vec 0\}\big)/\pm 1$ and the action of $t\in\mathbb R$ is given by $[\vec v]\to [e^t\vec v]$. Notice
that there is a projection map $\pi:T^1\mathbb H^2\to \partial^2\mathbb H^2$ (where $\partial^2\mathbb H^2=\{(x,y)\ :\ x,y\in\partial\mathbb H^2,\ x\ne y\}$) and
all the vectors tangent to the geodesic joining $x$ to $y$ are taken to $(x,y)$. Then, $\widetilde{{\bf U}_L}=\pi_*\widehat L$ is
a principal $\mathbb R$-bundle over $\partial^2\mathbb H^2$, so admits a natural geodesic flow. The group $\Gamma$ acts on  $\widetilde{{\bf U}_L}$
with quotient  ${\bf U}_L$ (see \cite[Prop. 2.4]{BCLS2} for more details). If $\gamma\in\Gamma$, then the closed orbit of ${\bf U}_L$ associated to $[\gamma]$
has period
$$-\log\frac{||\phi_{\ell(\gamma)}(\vec v)||}{||\vec v||}$$ 
for any vector $\vec v\in T^1\mathbb H^2$ tangent to the axis of $\gamma$. For example, if $L=L_1$, this period is the spectral radius 
$\omega_1(\rho(\gamma))=\log\lambda_1(\gamma)$ of $\rho(\gamma)$, while if $L=L_i\otimes L_j^*$, the period is given by 
$\alpha_{i,j}(\rho(\gamma))=\log\frac{\lambda_i(\rho(\gamma))}{\lambda_j(\rho(\gamma))}$. It is not difficult to write down an explicit
H\"older orbit equivalence between $T^1\mathbb H^2$ and ${\bf U}_L$ in general (see \cite[Prop. 4.2]{BCLS} for details).

\medskip

For more general $\phi\in\mathcal B^+(\rho)$, Sambarino \cite{sambarino-quantitative} makes use of the Iwasawa cocycle. Quint \cite{quint-ps} defines
the {\em Iwasawa cocycle} $B:\mathsf{PSL}(d,\mathbb R)\times\mathcal F_d\to \mathfrak{a}$ in terms of the Iwasawa decomposition. Specifically, if $F\in\mathcal F_d$,
then there exists $K\in\mathsf{PO}(d)$, so that $F=K(F_0)$ where $F_0$ is the flag determined by the standard basis and if $A\in\mathsf{SL}(d,\mathbb R)$, then
$B(A,F)$ is the unique element of $\mathfrak{a}$ satisfying
$$AK=Le^{B(A,F)} U$$
for some $L\in\mathsf{PO}(d)$ and upper triangular, unipotent element $U$. More geometrically, if $\vec v_1\in F^{(1)}$ is non-trivial, then 
$$\omega_1(B(A,F))=\log\frac{||A(\vec v_1)||}{||\vec v_1||}$$
and  if $\vec v_k$ is a non-trivial vector in $E^k(F^{(k)})\in\mathbb P(E^k\mathbb R^d)$ is non-trivial, where $E^k$ is the $k^{\rm th}$ exterior power, then
$$\omega_k(B(A,F))=\log\frac{||E^kA(\vec v_k)||}{||\vec v_k||}.$$
Notice that if $A$ is loxodromic and $F_A$ is the attracting flag of $A$, then $B(A,F_A)=\nu(A)$.
The Iwasawa cocycle satisfies the cocycle relation, $B(CD,F)=B(C,D(F))+B(D,F)$ for all $C,D\in\mathsf{PSL}(d,\mathbb R)$.

Given $\rho\in\mathcal H_d(S)$ and $\phi\in\mathcal B^+(\rho)$,  Sambarino defines the  H\"older cocycle
$$\beta_\rho^\phi:\Gamma\times\partial\mathbb H^2\to\mathfrak{a}\qquad\text{given by}\qquad \beta_\rho^\phi(\gamma,x)=\phi(B(\rho(\gamma),\xi_\rho(x))).$$
The period of a H\"older cocycle $\beta$ associated to $\gamma$ is  given by $\beta(\gamma,\gamma^+)$, so 
the period of $\beta_\rho^\phi$ associated to $\gamma$ is $\phi(\nu(\rho(\gamma))$. 
Using the theory of H\"older cocycle developed by Ledrappier \cite{ledrappier}, Sambarino shows that the action of $\Gamma$ on
$\partial^2\mathbb H^2\times\mathbb R$ defined by $\gamma(x,y,t)=(\gamma(x),\gamma(y), t+\beta_\rho^\phi(\gamma,y))$ is properly discontinuous and
cocompact. One may then define the quotient flow space ${\bf U}^\phi(\rho)=\partial^2\mathbb H^2\times\mathbb R/\Gamma$.

\medskip

Rafael Potrie and Andres Sambarino later  showed that simple root entropy is constant on the Hitchin component and used this to 
establish a remarkable entropy rigidity theorem for Hitchin representations.

\begin{theorem}{\rm (Potrie-Sambarino \cite{potrie-sambarino})}
If $\rho\in \mathcal H_d(S)$ and $1\le i\le d-1$, then $h^{\alpha_i}(\rho)=1$. Moreover, if $\phi=\sum c_i \alpha_i\in\Delta$, then
$$h^\rho(\phi)\le \frac{1}{c_1+\ldots+c_d}$$
and if $c_i>0$ for all $i$, then equality holds if and only if $\rho$ is $d$-Fuchsian.
\end{theorem}

\medskip\noindent
{\bf Remark:} Sambarino \cite[Cor. 7.15]{sambarino-quantitative} also obtains analogous results for the growth rate of translation length on the symmetric space.
In subsequent work, Sambarino \cite{sambarino-orbital} establishes a mixing property for the Weyl chamber flow which allowed him to establish equidistribution
results for ${\bf U}^\phi(\rho)$, see also Chow-Sarkar \cite{chow-sarkar}. For more recent developments, see Burger-Landesberg-Lee-Oh \cite{BLLO},
Carvajales \cite{carvajales-count, carvajales-growth},  Edwards-Lee-Oh \cite{ELO},
Landesberg-Lee-Lindenstraus-Oh \cite{LLLO}, Lee-Oh \cite{lee-oh, lee-oh2},
Pozzetti-Sambarino-Wienhard \cite{PSW,PSW2}, and Sambarino \cite{sambarino-dichotomy}.

\section{Pressure metrics for the Hitchin component}

In the 1970's, Bill Thurston proposed that one could construct a new Riemannian metric on Teichm\"uller space,
by considering the ``Hessian of the length of a random geodesic.'' Scott Wolpert  \cite{wolpert} (see also Fathi-Flaminio \cite{fathi-flaminio})
proved that Thurston's metric was a scalar multiple of the classical Weil-Petersson metric, which is defined using quadratic differentials
and Beltrami differentials. Bonahon \cite{bonahon} later re-interpreted Thurston's metric in terms of geodesic currents. McMullen \cite{mcmullen-wp}
showed that one could use the Thermodynamic formalism to construct a pressure form on the space $\mathcal H_0(T^1S)$ of all pressure zero H\"older functions on $T^1S$,
embed Teichm\"uller space in $\mathcal H_0(T^1S)$ and obtain Thurston's metric as the pullback of the pressure form. Bridgeman \cite{bridgeman-wp} extended 
McMullen's analysis to quasifuchsian space, obtaining a path metric which is an analytic Riemannian metric away from the Fuchsian locus. By construction, it
is mapping class group invariant and agrees with the Weil-Petersson metric (up to scalar multiplication) on the Fuchsian locus.

Martin Bridgeman, Dick Canary, Fran\c cois Labourie and Andres Sambarino \cite{BCLS,BCLS2} showed that one can use McMullen's procedure
to produce analytic pressure forms on the Hitchin component associated to any linear functional in $\Delta$. One key technical ingredient here is
to show that the limit map varies analytically over the space of Hitchin representations, see \cite[Theorem 6.1]{BCLS}. Since the Busemann cocycle $\beta_\rho^\phi$ is defined in
terms of the limit map, the Thermodynamic Formalism then implies that most natural dynamical quantities vary
analytically. In particular, the entropy varies analytically over the Hitchin component (see \cite[Thm. 1.3]{BCLS} and Pollicott-Sharp \cite[Thm.1.3]{pollicott-sharp-entropy}).
The foundational texts of Thermodynamic Formalism
are books by Bowen \cite{bowen-book}, Parry-Pollicott \cite{parry-pollicott} and Ruelle \cite{ruelle-book}. A fuller description of the use of Thermodynamic Formalism
to construct pressure metrics is given in the survey article by Bridgeman, Canary and Sambarino \cite{BCS}.

\medskip

Given two Hitchin representation $\rho:\Gamma\to\mathsf{PSL}(d,\mathbb R)$ and
$\eta:\Gamma\to\mathsf{PSL}(d,\mathbb R)$ and $\phi\in\Delta$, we define their {\em pressure intersection}
$$I^\phi(\rho,\eta)=\lim_{T\to\infty} \frac{1}{\#(R^\phi_T(\rho))}\sum_{[\gamma]\in R_T(\rho)}
\frac{\ell^\phi(\eta(\gamma))}{\ell^\phi(\rho(\gamma))}$$
which one may think of as the $\phi$-length (in $\eta$) of a random geodesic (with respect to $\phi$-length in $\rho$).
One then considers  the {\em renormalized pressure intersection} given by
$$J^\phi(\rho,\eta)=\frac{h^\phi(\eta)}{h^\phi(\rho)} I^\phi(\rho,\eta)$$
(Thurston and McMullen did not need to renormalize the pressure intersection since there is a single projective class of linear functionals and each entropy is constant on
the Teichm\"uller space of a closed surface.)
The functions $I^\phi$ and $J^\phi$ are analytic on the Hitchin component (see \cite[Thm. 1.3]{BCLS}) and $J^\phi$ achieves its global minimum of 1 along the diagonal
(see \cite[Cor. 8.2]{BCLS}). (The results in \cite{BCLS}  referenced in the last sentence are stated only for $\phi=\omega_1$, but the proofs easily generalize to
all linear functionals in $\Delta$, see \cite{BCKM2}.)

\begin{theorem}{\rm (Bridgeman-Canary-Labourie-Sambarino \cite[Cor. 6.2]{BCLS})}
If $S$ is a closed orientable surface, $d\ge 3$ and $\phi\in\Delta$, then $J^\phi$ is an analytic function on $\mathcal H_d(S)\times\mathcal H_d(S)$.
Moreover, if $\rho,\eta\in\mathcal H_d(S)$, then
$$J^\phi(\rho,\rho)=1\qquad\text{and}\qquad J^\phi(\rho,\eta)\ge 1.$$
\end{theorem}

If $\phi\in\Delta$, then one may define the {$\phi$-Pressure form} on $\mathcal H_d(S)$ by
$$\mathbb P^\phi|_{T_{\rho}\mathcal H_d(S)}=\mathrm{Hess}\ J^\phi(\rho,\cdot).$$
Since $J^\phi$ is analytic, $\mathbb P^\phi$ is analytic, and since $J^\phi$ achieves its minimum along
the diagonal, $\mathbb P^\phi$ is non-negative at every point. 
Notice that, by construction, $\mathbb P^\phi$ is mapping class group invariant and, by Wolpert's result \cite{wolpert},
agrees with (a scalar multiple of) the Weil-Petersson metric on the Fuchsian locus.
The most difficult portion of the analysis
then involves determining if $\mathbb P^\phi$ is non-degenerate, and hence gives rise to an analytic
Riemannian metric. 
(Mark Pollicott and Richard Sharp \cite{pollicott-sharp-pressure} provide an alternate formulation of the pressure form $\mathbb P^\phi$ when $\phi=\omega_1$.)

It is important to notice that  the $\phi$-pressure form will not be non-degenerate  for all $\phi$. For example,
since $\omega_H$ is invariant under the contragredient involution (which takes $A$ to $(A^T)^{-1}$), it is easy to see that $\mathbb P^{\omega_H}$
will be degenerate on the self-dual locus (i.e. the fixed point set of the contragredient involution) in $\mathcal H_d(S)$ (see \cite[Lem. 5.22]{BCS}),
which always includes the Fuchsian locus. The same analysis applies to any linear functional which is invariant under the contragredient involution.
For example, $\mathbb P^{\alpha_n}$ is degenerate on the self-dual locus of $\mathcal H_{2n}(S)$, see \cite[Prop. 8.1]{BCLS2}.
Similarly, the pressure metric on quasifuchsian space is degenerate  along the Fuchsian locus, which is the fixed point set of the involution
of quasifuchisan space induced by complex conjugation, see \cite{bridgeman-wp}.

However, in the case that $\phi$ is either the first fundamental weight $\omega_1$ or  first simple root $\alpha_1$,
$\mathbb P^\phi$ is non-degenerate. One hopes that Potrie and Sambarino's result that simple root entropy is constant
on the Hitchin component, will make the simple root pressure metric more tractable to study.
No other cases are fully understood at this point.

\begin{theorem}{\rm (\cite[Corollary 1.6]{BCLS} and \cite[Thm. 1.6]{BCLS2})}
If $S$ is a closed orientable surface and $d\ge 3$, then the pressure forms $\mathbb P^{\omega_1}$ and $\mathbb P^{\alpha_1}$ 
are analytic Riemannian metrics on the Hitchin component $\mathcal H_d(S)$ which are
invariant under the action of the mapping class group $\mathrm{Mod}(S)$.
Moreover, the restrictions of both metrics to the Fuchsian locus are scalar multiples
of the Weil-Petersson metric.
\end{theorem}

The following elementary conjecture illustrates how little is known about the pressure metric. See \cite[Sec. 7]{BCS} for a further discussion
of questions about the pressure metric, all of which remain wide open. However, since the time that survey was written, Fran\c cois Labourie
and Richard Wentworth \cite{labourie-wentworth} and Xian Dai \cite{dai} have made significant progress in describing the pressure metric
at the Fuchsian locus.

\medskip\noindent
{\bf Conjecture:} {\em
There exists a sequence of points in $\mathcal H_d(S)$ whose ${\mathbb P}^{\omega_1}$-distance to the Fuchsian locus
diverges to $\infty$.}

\medskip

One natural place to start would be to study Hitchin components of triangle groups when $d=3$, which are often one-dimensional, 
see Choi-Goldman \cite{choi-goldman-orbifold}.  (Nie \cite{nie-simplicial} explicitly computes the deformation spaces in some cases.) 
All the results discussed so far go through immediately for Hitchin components of
cocompact triangle groups (see, for example, Alessandrini-Lee-Schaffhauser \cite{ALS}).
It is still unknown whether or not these one-dimensional Hitchin components have finite diameter.
One would also like to investigate this conjecuture and all the questions in \cite{BCS} for the simple root pressure metric $\mathbb P^{\alpha_1}$.

\medskip

When $\phi=\omega_H$ and $d=3$,  one can
give a complete analysis of the degeneracy of $\mathbb P^{\omega_H}$. We say that a non-zero vector $\vec v\in T_\rho\mathcal H_d(S)$ is {\em self-dual}
if $dC(\vec v)=-\vec v$ where $C:\mathcal H_d(S)\to\mathcal H_d(S)$ is the contragredient involution, i.e. $C(\rho)\in\mathcal H_d(S)$ is given by
$C(\rho)(\gamma)=\rho(\gamma^{-1})^T$. Notice that self-dual vectors are based
at the self-dual locus and that when $d=3$ the self-dual locus is exactly the Fuchsian locus. 

\begin{theorem}
{\rm (Bridgeman-Canary-Sambarino \cite[Section 5.8]{BCS}, Bray-Canary-Kao-Martone \cite{BCKM2})}
If $S$ is a closed orientable surface and \hbox{$\vec v\in T\mathcal H_3(S)$} is non-zero, then
$\mathbb P^{\omega_H}(\vec v,\vec v)=0$ if and only if $\vec v$ is a self-dual vector. Therefore,
the pressure form $\mathbb P^{\omega_H}$ gives rise to  a mapping class group invariant path metric
which is an analytic Riemannian metric away from the Fuchsian locus  and agrees with (a scalar multiple of)
the Weil-Petersson metric on the Fuchsian locus.
\end{theorem}

\medskip\noindent
{\bf Remarks:} (1) Qiongling Li \cite{li-metric} produced another mapping class group invariant Riemannian
metric on $\mathcal H_3(S)$, which she calls the Loftin metric. The Loftin metric also restricts to a scalar multiple
of the Weil-Petersson metric on the Fuchsian locus. One expects that her metric differs from our pressure metrics,
but that is unknown so far.

Inkang Kim and Genkai Zhang \cite{kim-zhang} constructed a mapping class group invariant K\"ahler metric on $\mathcal H_3(S)$ in which the Fuchsian locus
is a totally geodesic complex submanifold whose intrinsic metric agrees with the Weil-Petersson metric, see also Labourie \cite[Cor. 1.3.2]{labourie-cyclic}.
The relationship of this metric to the pressure metrics and Li's metric is also not understood.

(2) Marc Burger \cite{burger} was the first one to consider the pressure intersection, in
the context of convex cocompact rank one  representations. His work was motivated by rigidity result for pairs of
Fuchsian representations due to Chris Bishop and  Tim Steger \cite{bishop-steger}. The pressure intersection
can also be interpreted in terms of Gerhard Knieper's geodesic stretch \cite{knieper}, see the discussion in Schapira-Tapie \cite{schapira-tapie}.

(3) Bridgeman, Canary, Labourie and Sambarino \cite{BCLS} more generally define a pressure
form $\mathbb P^{\omega_1}$ associated to the first fundamental weight at smooth points of
deformation spaces of projective Anosov representations into $\mathsf{SL}(d,\mathbb R)$ which is
non-degenerate at all ``generic'' representations.

\section{Geodesic currents and collar lemmas for Hitchin representations}

In this section, we discuss some of the work which further explores the analogy between the Hitchin component
and Teichm\"uller space. The choice of topics reflects our personal tastes and the focus of this article.

Francis Bonahon \cite{bonahon} exhibited a geodesic current $\mu_\rho$, known as the Liouville current,
associated to a Fuchsian representation $\rho:\pi_1(S)\to \mathsf{PSL}(2,\mathbb R)$ such that if
$\gamma\in\pi_1(S)-\{id\}$, then
$$i(\mu_\rho,\gamma)=\ell(\rho(\gamma))$$
(where $\ell(\rho(\gamma))$ is the hyperbolic
translation length of $\rho(\gamma)$). He used his theory of Liouville currents to reinterpret both Thurston's compactification of Teichm\"uller
space and Thurston's definition of the Weil-Petersson metric. Recall that a geodesic current on $S$ is
a locally finite, $\pi_1(S)$-invariant Radon measure on $\partial\pi_1(S)\times \partial \pi_1(S)- \Delta$
(where $\Delta$ is the diagonal). Moreover, currents associated to  (weighted) closed curves are dense in the space
$\mathcal C(S)$ of geodesic currents and  the intersection function 
$i:\mathcal C(S)\times\mathcal C(S)\to \mathbb R$
agrees with geometric intersection number on pairs  of closed curves.

If $\rho\in\mathcal H_d(S)$, Fran\c cois Labourie \cite{labourie-cross} (see also \cite{BCLS2} and Martone-Zhang \cite{martone-zhang}) exhibit a geodesic current
$\mu_\rho^H$, again called the {\em Liouville current}, so that 
$$i(\mu_\rho^H,\gamma)=\ell^{\omega_H}(\rho(\gamma))=\log\left(\frac{\lambda_1(\rho(\gamma))}{\lambda_d(\rho(\gamma))}\right)$$
if $\gamma \in\pi_1(S)-\{id\}$.
Bridgeman, Canary, Labourie, and Sambarino \cite{BCLS2} define the
{\em Liouville volume} 
$${\rm vol}_L(\rho)=i(\mu_\rho^H,\mu_\rho^H)$$
of a Hitchin representation. 
They also show that $\mu_\rho^H$ is a multiple
of the Bowen-Margulis current for the simple root flow $U^{\alpha_1}_\rho$ (see \cite[Thm. 1.3]{BCLS2}) and use this fact and
a result of Nicolas Tholozan \cite{tholozan-dominate} to establish a rigidity result when $d=3$.
It is natural to ask whether a similar rigidity holds in higher dimensions. 

\begin{theorem}{\rm (Bridgeman-Canary-Labourie-Sambarino \cite[Cor. 1.5]{BCLS2})}
If $S$ is a closed orientable surface and $\rho\in \mathcal H_3(S)$, then
${\rm vol}_L(\rho)\ge 4\pi^2|\chi(S)|$ with equality if and only if $\rho$ lies in the Fuchsian locus.
\end{theorem}

If $\mu$ is a geodesic current on $S$ and $\bf U$ is a geodesic flow orbit equivalent to $T^1S$, then one  may use the Hopf parametrization of $T^1\mathbb H^2$ to obtain a 
(possibly degenerate) volume form $\mu\otimes dt$  on
${\bf U}$ by considering the local product of $\mu$ and the element $dt$ of path length, see \cite{BCLS2} for details.
Bridgeman, Canary, Labourie and Sambarino show that one may re-interpret $I^{\alpha_1}=J^{\alpha_1}$ in terms of the Liouville current,
see \cite[Sec. 6]{BCLS2}. Specifically,
$$I^{\alpha_1}(\rho,\eta)=\frac{\int_{{\bf U}^{\alpha_1(\eta)}} \omega_\rho^H\otimes dt}{\int_{{\bf U}^{\alpha_1(\rho)}}\omega_\rho^H\otimes dt}.$$

More generally, Giuseppe Martone and Tengren Zhang \cite{martone-zhang} produce, for each $k\in\{1,\ldots,d-1\}$ and Hitchin representation $\rho:\pi_1(S)\to\mathsf{PSL}(d,\mathbb R)$, 
a geodesic current $\mu_\rho^k$ so that if $\gamma$ is a closed curve on $S$, then 
$$i(\mu_\rho^k,\gamma)=\ell^{\omega_k+\omega_{d-k}}(\rho(\gamma)).$$
In their construction, $\mu_\rho^1$ is exactly the symmetrization of $\mu_\rho^H$ (i.e the geodesic current with the same periods which is invariant under
the involution exchanging the first and second coordinate). Martone and Zhang also introduce the more general class of positively ratioed representations
and produce geodesic currents with analogous properties in this more general setting. 

Martone and Zhang use their theory to get a deep understanding of how the Hilbert length entropy of a sequence of Hitchin representations can
converge to 0. One particularly easy stated consequence of their work is the following relationship between the entropy and the systole of a Hitchin
representation. If $\phi\in\mathcal B^+(\rho)$, we can define the {\em $\phi$-systole} of $\rho$ to be 
$$\mathrm{sys}^\phi(\rho)=\min\big\{\ell^\phi(\gamma)\ |\ \gamma\in\Gamma\setminus\{id\}\big\}.$$

\begin{theorem}{\rm (Martone-Zhang \cite[Cor. 7.6]{martone-zhang})}
Let $S$ be a closed orientable surface.
Given $d\ge 3$ and $k\in\{1,\ldots,d-1\}$, there exists $L=L(S,d,k)>0$ so that if $\rho\in\mathcal H_d(S)$, then
$$h^{\omega_k+\omega_{d-k}}(\rho)\mathrm{sys}^{\omega_k+\omega_{d-k}}(\rho)\le L.$$
In particular, $h^{\omega_H}(\rho)\mathrm{sys}^{\omega_H}(\rho)\le L.$
\end{theorem}

Tengren Zhang \cite{zhang1,zhang2}  (and Xin Nie \cite{nie-simplicial,nie-entropy} when $d=3$) produce sequences 
$\{\rho_n\}$ in $\mathcal H_d(S)$ of representations  so that $\mathrm{sys}^{\omega_H}(\rho_n)\to\infty$ and $h^{\omega_H}(\rho_n)\to 0$.
%, where
%$\omega_H=\omega_1+\omega_{d-1}$.
If $d=3$,  $\omega^H(\rho_n(\gamma))=\ell^{\alpha_1}(\rho_n(\gamma))+\ell^{\alpha_1}(\rho_n(\gamma^{-1})$, so  
if $\gamma$ is any element of $\pi_1(S)-\{1\}$, then  
$\max\{\ell^{\alpha_1}(\rho_n(\gamma)),\ell^{\alpha_1}(\rho_n(\gamma^{-1}))\}\to \infty$. 
One hopes that such sequences would have their distance to the Fuchsian locus diverge to infinity.

One might naively think that Zhang's sequences would have $h^{\alpha_1}(\rho_n)\to 0$.
However, we know, from Potrie-Sambarino \cite{potrie-sambarino},
that $h^{\alpha_1}(\rho_n)=1$ for all $n$. One possible explanation for this, in the simple case where $d=3$,
would be that for many elements in $\pi_1(S)$, the middle eigenvalue  $\lambda_2(\rho_n(\gamma))$  remains ``near'' to either 
$\lambda_1(\rho_n(\gamma))$ or $\lambda_3(\rho_n(\gamma))$,
so that one of $\{\ell^{\alpha_1}(\rho_n(\gamma)),\ell^{\alpha_1}(\rho_n)(\gamma^{-1})\}$ is growing quickly
while the other remains moderate or even bounded. This phenomenon was observed in explicit computations done
by Martin Bridgeman and the author for Hitchin components of certain triangle groups.
This suggests the surprising possibility that the simple root systole is bounded above on $\mathcal H_d(S)$.

\medskip\noindent
{\bf Question:} {\em
Given a closed surface $S$ and $d$,  does there exist $L>0$ so that 
${\rm sys}^{\alpha_1}(\rho)\le L$ for all $\rho \in\mathcal H_d(S)$?}

\medskip

Another deep analogy with the traditional theory of Fuchsian groups was established when
Gye-Seon Lee and Tengren Zhang \cite{lee-zhang}  proved an analogue of the collar lemma for Fuchsian groups for Hitchin representations.
We will not state the precise version of their results,
but we recall the following consequences which hold for all Hitchin representations.

\begin{theorem}{\rm (Lee-Zhang \cite[Cor. 1.2]{lee-zhang})}
If $S$ is a closed, orientable surface of genus at least two, $\alpha$ and $\gamma$ are homotopically non-trivial closed
curves on $S$, and $\rho\in\mathcal H_d(S)$, then
\begin{enumerate}
\item
if $i(\alpha,\gamma)\ne 0$, then 
$$\left(e^{\ell^{\omega_H}(\rho(\alpha))}-1\right)\left(e^{\ell^{\omega_H}(\rho(\alpha))}-1\right)>1,$$
\item
and if $\alpha$ is not simple, then $\ell^{\omega_H}(\rho(\alpha))\ge\log(2)$.
\end{enumerate}
\end{theorem}

One consequence of their work is that sufficiently complicated curve systems determine proper multi-length functions.

\begin{corollary}{\rm (Lee-Zhang \cite[Cor. 1.4]{lee-zhang})}
If $S$ is a closed orientable surface of genus at least two and $\Gamma=\{\gamma_1,\ldots,\gamma_k\}$ is a collection of homotopically
non-trivial closed curves which contains a pants decomposition of $S$ and so that  if $\alpha$ is any homotopically non-trivial curve $\alpha$ on $S$,
then $i(\alpha,\Gamma)\ne 0$, then the map
$$L_\Gamma:\mathcal H_d(S)\to \mathbb R^k\text{ given by } L_\Gamma(\rho)=\left(\ell^{\omega_H}(\rho(\gamma_i))\right)_{i=1}^k$$
is proper.
\end{corollary}

Marc Burger and Beatrice Pozzetti \cite{burger-pozzetti} subsequently established collar lemmas for maximal representations into $\mathsf{Sp}(2n,\mathbb R)$ and
Jonas Beyrer and Beatrice  Pozzetti \cite{beyrer-pozzetti-collar} for partially hyperconvex representations of surface groups.

\medskip

We briefly mention some other important work which explore analogies between Hitchin components and Teichm\"uller spaces.

\begin{enumerate}
\item
Richard Skora \cite{skora} showed that one may describe Thurston's compactification of Teichm\"uller space in terms of actions of $\pi_1(S)$ on $\mathbb R$-trees.
Anne Parreau \cite{parreau} constructed a similar compactification of the Hitchin componen by actions of
$\pi_1(S)$ on  $\mathbb R$-buildings. (Her work applies to compactification of  much more general character varieties.) Marc Burger, 
Alessandra Iozzi,  Anne Parreau and  Beatrice Pozzetti \cite{BIPP-compact}
analyzed the Parreau compactification and showed that there is a non-empty open domain of discontinuity for the action of the mapping class group on the boundary.
In subsequent work, they study the real spectrum compactification of the Hitchin component which admits a continuous, surjective map to the Parreau compactification which
is mapping class group invariant, see their survey article \cite{BIPP-survey}. Their work depends crucially on the theory of geodesic currents.
\item
Zhe Sun, Anna Wienhard and Tengren Zhang \cite{SWZ,SZ,WZ} studied Goldman's symplectic form on the Hitchin component. They prove
an analogue of Wolpert's magic formula in this setting and construct Darboux coordinates for the symplectic structure. They also construct a half-dimensional
space of Hamiltonian flows which generalize the twist flows on Teichm\"uller space.
\item
Nigel Hitchin \cite{hitchin} offers a paremetrization of $\mathcal H_d(S)$ as $\prod_{d=2}^k Q^k(X)$ where 
$Q^k(X)$ is the space of holomorphic $k$-differentials on a Riemann surface $X$ homeomorphic to $S$. Unfortunately, this parameterization is
not invariant with respect to the mapping class group. Fran\c cois Labourie \cite{labourie-cubic} and John Loftin \cite{loftin-cubic} showed that there is a homeomorphism from
$\mathcal H_3(S)$ to the bundle of cubic holomorphic differentials over $\mathcal T(S)$ which is equivariant with respect to the action of the mapping class group.
Recent work of Vlad Markovic \cite{markovic-nonunique} indicates that  Labourie's approach in \cite{labourie-cyclic} will not generalize to produce similar  mapping class group invariant parametrizations
when $d>3$.
\item
Bridgeman, Pozzetti, Sambarino and Wienhard \cite[Cor. A]{BPSW} showed that $\mathcal H_d(S)$ is an isolated minimum for the $\alpha_1$-entropy functional for
the space of quasi-Hitchin representations into $\mathsf{PSL}(d,\mathbb C)$. (A representation $\rho:\pi_1(S)\to\mathsf{PSL}(d,\mathbb C)$ is quasi-Hitchin
if it is Borel Anosov and can  be deformed to a Hitchin representation through Borel Anosov representations.) Moreover, they describe the Hessian of the $\alpha_1$-entropy functional
at the Hitchin locus. Their results generalize work of Bowen \cite{bowen-qf}, Bridgeman \cite{bridgeman-wp} and McMullen \cite{mcmullen-wp} in the quasifuchsian setting.
\item
Fran\c cois Labourie and Greg McShane \cite{labourie-mcshane} proved an analogue of McShane's identity for Hitchin representations, while 
Nick Vlamis and Andrew Yarmola \cite{vlamis-yarmola} proved an analogue of the Basmajian identity for Hitchin representations. Yi Huang and Zhe Sun \cite{huang-sun}
proved versions of McShane's identity for holonomy maps of finite area convex projective surfaces (and for positive representations of Fuchsian lattices).
\item
Richard Schwartz and Richard Sharp \cite{schwartz-sharp} proved an explicit correlation result for lengths on hyperbolic surfaces. Specifically, they prove that 
given $\rho_1,\rho_2\in\mathcal T(S)=\mathcal H_2(S)$ and $\epsilon>0$, there exist $C,M>0$ so that
 $$\#\Big\{[\gamma]\in [\Gamma]\ :\ \ell(\rho_1(\gamma))\in (x,x+\epsilon)\  \text{ and }\ \ell(\rho_2(\gamma))\in (x,x+\epsilon)\Big\}\sim C\frac{e^{Mx}}{x^{3/2}}.$$
Xian Dai and Giuseppe Martone \cite[Thm. 1.7]{dai-martone} generalize this result to arbitrary linear functions in $\Delta$, by showing that given $\rho_1,\rho_2\in\mathcal H_d(S)$, $\phi\in\Delta$
and $\epsilon>0$, there exist $C,M>0$ so that
{\footnotesize
$$\#\Big\{[\gamma]\in [\Gamma]\ :\ h^\phi(\rho_1)\ell^\phi(\rho_1(\gamma))\in (x,x+h^\phi(\rho_1)\epsilon)\  \text{ and }\ h^\phi(\rho_2)\ell^\phi(\rho_2(\gamma))\in (x,x+h^\phi(\rho_2)\epsilon)\Big\}\sim C\frac{e^{Mx}}{x^{3/2}}.$$}
Schwartz and Sharp \cite{schwartz-sharp}  asked whether or not $M$ can be arbitrarily close to 0. Dai and Martone  \cite[Thm. 1.3]{dai-martone} show that this can happen even in the Fuchsian setting.
\item
It is a classical result, that one can find finitely many simple closed curves whose lengths determine a point in Teichm\"uller space. Ursula Hamenstadt \cite{hamenstadt} 
and Paul Schmutz \cite{schmutz} showed
that $6g-5$ curves suffice, but that no set of $6g-6$ curves completely determine a point in Teichm\"uller space. 
Bridgeman, Canary, Labourie and Sambarino \cite[Thm 1.2]{BCLS} show that the $\omega_1$-lengths of all curves on $S$ determine a point in $\mathcal H_d(S)$.
Bridgeman, Canary and Labourie \cite{BCL} showed that the $\omega_1$-length of
all simple closed curves on $S$ determine a points in $\mathcal H_d(S)$ if $S$ has genus at least 3. It would be interesting to know whether or not finitely many curves suffice.
(Sourav Ghosh \cite{ghosh} recently showed that there are finitely many elements of $\pi_1(S)$ whose full Jordan projections determine a point in $\mathcal H_d(S)$.)
\item
Vladimir Fock and Alexander Thomas  \cite{fock-thomas} have defined the notion of a higher complex structure on a closed surface. They conjecture that the
space of higher complex structures of order $d$ on a surface $S$  is canonically isomorphic to  the Hitchin component $\mathcal H_d(S)$.

\end{enumerate}

\section{The Teichm\"uller theoretic viewpoint of Sullivan and Tholozan}

In August 2017, Nicolas Tholozan gave an inspirational talk \cite{tholozan-talk} describing how to use work of Dennis Sullivan
to give a Teichm\"uller-theoretic interpretation of the simple root pressure metric $\mathbb P^{\alpha_1}$. We will give a brief description of his work, which unfortunately
is not fully available yet. If you want more details I suggest you view his talk \cite{tholozan-talk} on YouTube, consult his lecture notes \cite{tholozan-notes} for
a mini-course given at the University of Michigan in December 2019 and/or read Sullivan's beautiful paper \cite{sullivan-universal}.

One may consider the unit tangent bundle $T^1\mathbb H^2$ of the hyperbolic plane as a foliated space where the leaves are the central stable leaves of the geodesic flow.
More prosaically, each leaf of the foliation consists of tangent vectors to all geodesic ending at a fixed point in $\partial\mathbb H^2$.
Each leaf is canonically identified with $\mathbb H^2$ (up to isometry) and hence admits a well-defined complex structure. One may thus view
$T^1\mathbb H^2$ as admitting a complex foliation. If $S=\mathbb H^2/\Gamma$, then $\Gamma$ acts as a group of holomorphic automorphisms of
this complex foliation, so $T^1S$ is also a complex foliated manifold.  

Sullivan \cite{sullivan-universal} developed a theory of the Teichm\"uller space
$\mathcal T(L)$ of complex structures on
a complex laminated space $L$ (i.e. spaces which admit a local product structure in which the horizontal leaves admit a complex structure).
Notice that a complex structure on a leaf induces a smooth structure on the leaf.
An element of $\mathcal T(L)$ may be viewed as a complex foliated laminated space with a homeomorphism to $L$ which
preserves the lamination structure and is smooth on each leaf (up to appropriate marked equivalence). Sullivan's Teichm\"uller theory involves generalizations of 
both quadratic differentials and Beltrami differentials, so he is able to construct both a Teichm\"uller metric and a Weil-Petersson metric on $\mathcal T(L)$.

Tholozan \cite{deroin-tholozan} constructs a continuous bijection $CF$ between the space of (conjugacy classes of) entropy one geodesic flows which
are H\"older orbit equivalent to $T^1S$ and the space $\mathcal T^h(T^1S)$ of elements of $\mathcal T(T^1S)$ where the homeomorphism to $T^1(S)$ is
transversely H\"older, see \cite[Thm. 0.4]{tholozan-notes}. 
Roughly, one maps the stable leafs of flow space to horocycles. The fact that the entropy is $1$ allows one to see that this gives an identification
of (the cover of) each leaf with $\mathbb H^2$. A result of Candel \cite{candel} then shows that this gives
rise to a foliated complex structure on $T^1S$.

From $CF$ one obtain an embedding
$$R:\mathcal H_d(S)\to\mathcal  T(T^1S)$$
where $R(\rho)=CF({\bf U}^{\alpha_1}(\rho))$.
Tholozan's main result is that $\mathbb P^{\alpha_1}$ is the pullback of a  ``Weil-Peterson metric''
on $\mathcal T(T^1S)$.

\begin{theorem}{\rm (Tholozan \cite{deroin-tholozan})}
If $S$ is a closed orientable hyperbolic surface of genus $g\ge 2$, then the simple root pressure metric $\mathbb P^{\alpha_1}$ is the pull-back, via
the embedding $R$, of a scalar multiple of a ``Weil-Petersson metric'' on $\mathcal T(T^1S)$.
\end{theorem}

This opens up the possibility of using classical Teichm\"uller theoretic techniques to study the, so far mysterious, simple root pressure metric.
One can also pull back the Teichm\"uller metric on $\mathcal T(T^1S)$ to obtain a metric $\mathbb Q^{\alpha_1}$ on $\mathcal H_d(S)$,
which one might call the {\em simple root Teichm\"uller metric}.
As in the classical setting,  the simple root Teichm\"uller metric should be less regular, but easier to control. One might first study the properties of the simple root Teichm\"uller metric,
and then study its relationship with the simple root pressure metric.

\medskip\noindent
{\bf Remark:} Labourie \cite{labourie-cross} also constructs a candidate for a highest Teichm\"uller space in which all Hitchin components embed.
He uses cross ratios to embed every Hitchin component $\mathcal H_d(S)$ in the space of (conjugacy classes of) representations of $\pi_1(S)$  into the space of
H\"older self-homeomorphisms of the space $J$ of 1-jets of real-valued functions on the circle. This embedding records the Hilbert length functional and is
closely related to Labourie's Liouville current. There is a relationship between Labourie's work and that of Tholozan coming from the fact that the Liouville
current is the Bowen-Margulis current of the (first) simple root flow.

\section{Hitchin representations of Fuchsian groups}

If $\Gamma$ is a Fuchsian group, i.e. a discrete subgroup of $\mathsf{PSL}(2,\mathbb R)$, we say that 
a representation $\rho:\Gamma\to\mathsf{PSL}(d,\mathbb R)$ is a Hitchin representation if there exists a $\rho$-equivariant positive map
$\xi:\Lambda(\Gamma)\to \mathcal F_d$, where $\Lambda(\Gamma)\subset\partial\mathbb H^2$ is the limit set of $\Gamma$.  
If $\Gamma$ is convex cocompact (i.e. finitely generated and without parabolic elements) and torsion-free, Hitchin representations of
$\Gamma$ were studied by Labourie and McShane \cite{labourie-mcshane}.
If $\mathbb H^2/\Gamma$ has finite volume and $d=3$, then
$\rho(\Gamma)$ preserves and acts properly discontinuously on a strictly convex domain $\Omega_\rho\subset\mathbb{RP}^2$ and
$\Omega_\rho/\rho(\Gamma)$ is a finite area real projective surface (or orbifold), see Choi-Goldman \cite{goldman-choi} or Marquis \cite{marquis-finite,marquis-deformation}.

Dick Canary, Tengren Zhang and Andy Zimmer \cite{CZZ} prove that Hitchin representations of finitely generated Fuchsian groups have many
of the same geometric properties as Hitchin representations of closed surface groups. (More generally, they study Anosov representations of
finitely generated Fuchsian groups.)

\begin{theorem} {\rm (Canary-Zhang-Zimmer \cite{CZZ})}
\label{cusped Hitchin properties}
If $\Gamma$ is a geometrically finite Fuchsian group, $b_0$ is a basepoint for $\mathbb H^2$  and  $\rho:\Gamma\to\mathsf{PSL}(d,\mathbb R)$ is a Hitchin representation, then
\begin{enumerate}
\item
There exists $B,b>0$ so so that if $\gamma\in\Gamma$, then
$$Be^{b\ell(\gamma)}\ge \frac{\lambda_k(\rho(\gamma))}{\lambda_{k+1}(\rho(\gamma))}\ge \frac{1}{B}e^{\frac{\ell(\gamma)}{b}}$$
for all $k\in\{1,\ldots,d-1\}$.
\item
If $\alpha\in\Gamma$ is parabolic, then $\rho(\alpha)$ is unipotent and its Jordan normal form has only one block.
\item
The orbit map $\tau_\rho:\Gamma(b_0)\to X_d(\mathbb R)$ given by $\tau_\rho(\gamma(b_0))=[\rho(\gamma)]$ is a quasi-isometric embedding.
\item
The limit map $\xi_\rho$ is H\"older.
\item
If $z\in\Lambda(\Gamma)$, then $\xi_\rho(z)$ varies analytically over the space of  Hitchin representations.
\end{enumerate}
\end{theorem}

Harry Bray, Dick Canary, Nyima Kao and Giuseppe Martone \cite{BCKM} developed analogues of the dynamical results of Sambarino
in the setting of general Hitchin representations of  a torsion-free, finitely generated Fuchsian group $\Gamma$. 
If $\Gamma$ is convex cocompact, then Sambarino's original theory applies. The key new difficulty in the presence of parabolic elements
is that one can no longer model the (recurrent portion of the) geodesic flow on $T^1X$, where $X=\mathbb H^2/\Gamma$ by a finite
Markov coding. However, Fran\c coise Dal'bo and Marc Peign\'e \cite{dalbo-peigne} (in the case where $X$ has infinite area) and Manuel Stadlbauer, Francois Ledrappier
and Omri Sarig  \cite{stadlbauer,ledrappier-sarig} (in the case where $X$ has finite area) have developed well-behaved countable Markov codings.
These codings are natural generalizations of the Bowen-Series \cite{bowen-series} (finite) Markov codings of convex cocompact Fuchsian groups.
Kao \cite{kao-pm} first used these codings to construct a pressure metric on Teichm\"uller space, and Bray, Canary and Kao \cite{BCK} used them to
construct pressure metrics on deformation spaces of cusped quasifuchsian groups.

We recall that a one-sided countable Markov shift $(\Sigma^+,\sigma)$ is determined by a countable alphabet $\mathcal A$ and a transition
matrix $T\in\{0,1\}^{\mathcal A\times\mathcal A}$. An element $x\in\Sigma^+ $ is a one-sided infinite string $x=(x_i)_{i\in\mathbb N}$ of letters in $\mathcal A$
so that $T(x_i,x_{i+1})=1$ for all $i\in\mathbb N$. The shift $\sigma$ simply removes the first letter and shifts every other letter one place to the left,
i.e. $\sigma (x)=(x_{i+1})_{i\in\mathbb N}$.  Let $\mathrm{Fix}^n$ denote the set of periodic words in $\Sigma^+$ with period $n$.

Given a torsion-free finitely generated group $\Gamma$, then the associated countable Markov shift $(\Sigma^+,\sigma)$ constructed by Dal'bo-Peign\'e or Stadlbauer-Ledrappier-Sarig
has the following crucial properties.
\begin{enumerate}
\item
There exists a  finite-to-one H\"older map $\omega:\Sigma^+\to\Lambda(\Gamma)$ which surjects onto the complement $\Lambda_c(\Gamma)$
of the set of fixed points of parabolic elements of $\Gamma$.
\item
There exists a map $G:\mathcal A\to \Gamma$ so that if $x\in\mathrm{Fix}^n$, then $\omega(x)$ is the attracting fixed point of $G(x_1)\cdots G(x_n)$.
\item
If $\gamma\in\Gamma$ is hyperbolic, then there exists $x\in\mathrm{Fix}^n$ (for some $n$) so that $\gamma$ is conjugate to $G(x_1)\cdots G(x_n)$. Moreover, $x$
is unique up to cyclic permutation.
\end{enumerate}
In addition, $\Sigma^+$
is well-behaved in the sense that it satisfies the assumptions needed to make use of the powerful Thermodynamic Formalism for countable Markov shifts developed by
Daniel Mauldin and Mariusz Urbanski \cite{MU} and Omri Sarig \cite{sarig-first,sarig-2009}.

Bray, Canary, Kao and Martone \cite{BCKM} show that given $\rho\in\mathcal H_d(\Gamma)$ there is a vector-valued function $\tau_\rho:\Sigma^+\to\mathfrak{a}$
which records all the spectral data of $\rho(\Gamma)$. Specifically, they define
$$\tau_\rho(x)=B\big(\rho(G(x_1)),\rho(G(x_1))^{-1}(\xi_\rho(\omega(x)))\big)$$
and prove that it has the following property.

\begin{theorem}{\rm (Bray-Canary-Kao-Martone \cite{BCKM})}
Suppose that $\Gamma$ is a torsion-free, finitely generated Fuchsian group, with associated one-sided Markov shift $(\Sigma^+,\sigma)$ and 
$\rho:\Gamma\to\mathsf{PSL}(d,\mathbb R)$ is a Hitchin representation. 
There is a locally H\"older continuous function $\tau_\rho:\Sigma^+\to \mathfrak{a}$ so that if $\phi\in\Delta$, $\tau_\rho^\phi=\phi\circ\tau_\rho$, and
$x\in\mathrm{Fix}^n$, then 
$$S_n\tau_\rho^\phi(x)=\sum_{i=0}^{n-1} \tau_\rho^\phi(\sigma^i(x))=\ell^\phi(G(x_1)\cdots G(x_n)).$$
\end{theorem}

Much as in the classical case one may define the $\phi$-topological entropy of a Hitchin representations by letting
$$R_T^\phi(\rho)=\{[\gamma]\in[\Gamma_{hyp}]\ :\ \ell^\phi(\rho(\gamma))\le T\} \qquad\text{and}\qquad  h^\phi(\rho)=\lim \frac{\log\#R_T^\phi(\rho)}{T}$$
where $[\Gamma_{hyp}]$ is the collection of conjugacy classes of  hyperbolic elements of $\Gamma$. The only difference here is that we omit
consideration of parabolic elements, which are not present in the case of closed surface groups. One may then also generalize the definitions of pressure
intersection and renormalized intersection. If $\rho,\eta\in\mathcal H_d(\Gamma)$ and $\phi\in\Delta$, then
$$I^\phi(\rho,\eta)=\lim_{T\to\infty} \frac{1}{\#(R^\phi_T(\rho))}\sum_{[\gamma]\in R_T(\rho)}
\frac{\ell^\phi(\eta(\gamma))}{\ell^\phi(\rho(\gamma))}$$
and
$$J^\phi(\rho,\eta)=\frac{h^\phi(\eta)}{h^\phi(\rho)} I^\phi(\rho,\eta).$$

Bray, Canary, Kao and Martone use the renewal theorem of Marc Kesseb\"ohmer and Sabrina Kombrink \cite{kess-komb} for countable Markov shifts to establish counting and equidistribution
results in the setting of countable Markov shifts in the spirit of the work of Steven Lalley \cite{lalley} in the setting of finite Markov shifts. In the case of Hitchin
representations, their counting result has the following form. (When $d=3$, then our results are a special case of more general results of
Feng Zhu  \cite{zhu-dynamics} when $\phi=\omega_H$.)

\begin{theorem}{\rm (Bray-Canary-Kao-Martone \cite{BCKM})}
Suppose that $\Gamma$ is a torsion-free, finitely generated Fuchsian group and 
$\rho:\Gamma\to\mathsf{PSL}(d,\mathbb R)$ is a Hitchin representation. 
If $\phi\in\mathcal B^+(\rho)$, then there exists $h^\phi(\rho)>0$ so that
$$\#R_T^\phi(\rho)\sim \frac{e^{t h^\phi(\rho)}}{t h^\phi(\rho)}.$$
\end{theorem}

Our equidistribution result  \cite[Cor. 1.6]{BCKM} expresses the geometrically defined pressure intersection function in
terms of equilibrium states, which allows one to employ the machinery of Thermodynamical Formalism to verify analyticity and construct
analytic pressure forms. A more precise statement would take us further into a discussion of Thermodynamic Formalism than time allows
for in a brief survey paper.

The results in \cite{BCKM} and \cite{CZZ} combine to show that entropies and pressure intersection vary analytically.

\begin{corollary} {\rm (\cite{BCKM2})}
\label{hitchin analyticity}
If $\Gamma$ is a finitely generated Fuchsian group and $\phi\in\Delta$, then
$h^\phi$ is an analytic function on $\mathcal H_d(\Gamma)$ and $I^\phi$ and $J^\phi$ are analytic functions
on $\mathcal H_d(\Gamma)\times\mathcal H_d(\Gamma)$. Moreover, $J^\phi(\rho,\rho)=1$ and
$J^\phi(\rho,\eta)\ge 1$ for all $\rho,\eta\in\mathcal H_d(\Gamma)$.
\end{corollary}

If $\phi\in\Delta$, then we can again define a $\phi$-pressure form $\mathbb P^{\phi}$ on $\mathcal H_d(\Gamma)$ by considering
the Hessian of the renormalized pressure intersection $J^\phi$.
Corollary \ref{hitchin analyticity} allows one to use the Thermodynamic Formalism for countable Markov shifts developed
by Mauldin-Urbanski and Sarig to produce pressure metrics on Hitchin components of general finitely generated, torsion-free Fuchsian groups.
(One can embed Hitchin components of finitely generated Fuchsian groups with torsion into Hitchin components of finitely generated, torsion-free
Fuchsian groups and thus obtain pressure metrics on them as well.) Recall that $\mathrm{Mod}(\Gamma)$ is the group of (isotopy classes of) orientation-preserving self-homeomorphisms of
$X=\mathbb H^2/\Gamma$.

\begin{theorem}{\rm (Bray-Canary-Kao-Martone \cite{BCKM2})}
If $\Gamma$ is a finitely generated torsion-free Fuchsian group and  $d\ge 3$, then the pressure forms $\mathbb P^{\omega_1}$ is
an analytic Riemannian metric on the Hitchin component $\mathcal H_d(\Gamma)$ which is
invariant under the action of the mapping class group $\mathrm{Mod}(\Gamma)$.
\end{theorem}

Canary, Zhang and Zimmer \cite{CZZ2} showed that if $\Gamma$ is a lattice (i.e. $\mathbb H^2/\Gamma$ has finite area), then
the simple root entropies are constant over $\mathcal H_d(\Gamma)$, which generalizes a result of Potrie-Sambarino \cite{potrie-sambarino}
from the cocompact case.

\begin{theorem}{\rm (Canary-Zhang-Zimmer \cite{CZZ2})}
If $\Gamma$ is a torsion-free Fuchsian lattice and  \hbox{$\rho\in\mathcal H_d(\Gamma)$},
then $h^{\alpha_k}(\rho)=1$ for all $k\in\{1,\ldots,d-1\}$.
\end{theorem}

This allows Bray, Canary, Kao and Martone to establish the non-degeneracy of the simple root pressure metric.

\begin{theorem}{\rm (Bray-Canary-Kao-Martone \cite{BCKM2})}
If $\Gamma$ is a torsion-free Fuchsian lattice and  $d\ge 3$, then the pressure form $\mathbb P^{\alpha_1}$ is
an analytic Riemannian metric on the Hitchin component $\mathcal H_d(\Gamma)$ which is
invariant under the action of the mapping class group $\mathrm{Mod}(\Gamma)$.
\end{theorem}

The analysis of the Hilbert length pressure metric was carried out for all finitely generated Fuchsian groups,
yielding the following general result.

\begin{theorem}
{\rm (Bray-Canary-Kao-Martone \cite{BCKM2})}
If $\Gamma$ is a finitely generated torsion-free Fuchsian group and 
$\vec v\in T_\rho\mathcal H_3(\Gamma)$ is non-zero, then
$\mathbb P^{\omega_H}(\vec v,\vec v)=0$ if and only if $\vec v$ is a self-dual vector. Therefore,
the pressure form $\mathbb P^{\omega_H}$ gives rise to  a mapping class group invariant path metric
which is an analytic Riemannian metric away from the Fuchsian locus.
\end{theorem}

\medskip

Canary, Zhang and Zimmer \cite{CZZ2} are also able to generalize the entropy rigidity theorem of Rafael Potrie and Andres Sambarino \cite{potrie-sambarino}.
(In the process they obtain a result for the Hausdorff dimension of $(1,1,2)$-hypertransverse groups which is a common
generalization of the results of Beatrice Pozzetti, Andres Sambarino and Anna Wienhard \cite{PSW} for $(1,1,2)$-hyperconvex Anosov representations
and those of Chris Bishop and Peter Jones \cite{bishop-jones} for discrete subgroups of $\mathsf{SO}(d,1)$.) We recall that Sambarino \cite{sambarino-positive},
generalizing earlier unpublished work of Olivier Guichard, showed that if $\rho:\Gamma\to\mathsf{PSL}(d,\mathbb R)$ is Hitchin, then either
$\rho(\Gamma)$ is Zariski dense or its Zariski closure is conjugate to either $\tau_d(\mathsf{PSL}(2,\mathbb R))$ (in which case $\rho$ is
$d$-Fuchsian), $\mathsf{G}_2$ (in which case $d=7$), $\mathsf{PSO}(n,n-1)$ (in which case $d=2n-1$) or $\mathsf{PSp}(2n,\mathbb R)$
(in which case $d=2n$).

\begin{theorem}{\rm (Canary-Zhang-Zimmer \cite{CZZ2})}
If $\Gamma$ is a finitely generated  Fuchsian group, \hbox{$\rho\in\mathcal H_d(\Gamma)$} and $\phi=\sum c_i\alpha_i\in\Delta$, then
$$h^\phi(\rho)\le \frac{1}{c_1+\ldots+c_{d-1}}.$$
Moreover, equality  occurs exactly when $\Gamma$ is a lattice and either
\begin{enumerate}
\item
$\rho$ is $d$-Fuchsian
\item
$\phi=c_k\alpha_k$ for some $k$.
\item
$d=2n$, the Zariski closure of $\rho(\Gamma)$ is conjugate into $\mathsf{PSp}(2n,\mathbb R)$ and $\phi=c_k\alpha_k+c_{d-k}\alpha_{d-k}$ for some $k$.
\item
$d=2n-1$, the Zariski closure of $\rho(\Gamma)$ is conjugate  into $\mathsf{PSO}(n,n-1)$ and  $\phi=c_k\alpha_k+c_{d-k}\alpha_{d-k}$ for some $k$.
\item
$d=7$, the Zariski closure of $\rho(\Gamma)$ is conjugate to $\mathsf{G}_2$,  and $\phi=c_1\alpha_k+c_3\alpha_3+c_4\alpha_4+c_6\alpha_6$.
\end{enumerate}
\end{theorem}

\medskip\noindent
{\bf Remarks:} (1) Hitchin representations of finitely generated groups are relatively dominated, in the sense developed by Feng Zhu \cite{zhu-reldom,zhu-reldom2},
and relatively Anosov from various viewpoints developed in the work of Misha Kapovich and Bernhard Leeb \cite{kapovich-leeb}. Many of the properties
detailed in Theorem \ref{cusped Hitchin properties} can be derived in their frameworks. Our main motivation for developing our viewpoint was to 
establish the analytic variation of the limit map, which was not yet available from either of the previous viewpoints.

(2) Our work with Bray, Kao and Martone, was partially inspired by the work of Barbara Schapira and Samuel Tapie \cite{schapira-tapie}. In particular,
their work develops the notion of an entropy gap at infinity for a geodesic flow on a negatively curved manifold, which we adapt in our setting of
H\"older potentials on ``well-behaved'' countable Markov shifts, see also Velozo \cite{velozo}. One may obtain related counting and equidistribution
results for cusped Hitchin representations using the theory in Paulin-Pollicott-Schapira \cite{PPS} and/or Schapira-Tapie \cite{schapira-tapie-counting}.

\section{The augmented Hitchin component}

In this section, we recall the theory of the augmented Teichm\"uller space of a closed surface from classical Teichm\"uller theory and describe
an analogous conjectural geometric picture of the augmented Hitchin component.  The Hitchin component $\mathcal H_3(S)$ can be viewed
as the space of (marked) real projective surfaces homeomorphic to $S$, so the analogy is easiest to discuss when $d=3$. Moreover, John Loftin
and Tengren Zhang \cite{loftin-zhang} have explored the topological picture when $d=3$. We will restrict our discussion to this case, although we hope that there is an
analogous picture for all $d$.

\medskip\noindent
{\bf Augmented Teichm\"uller space:}
The augmented Teichm\"uller space $\widehat{\mathcal T}(S)$ of a closed orientable surface $S$ of genus $g\ge 2$, is obtained from Teichm\"uller space by appending
all finite area hyperbolic surfaces obtained by pinching a collection of disjoint simple closed curves on $S$. It was introduced by Bill Abikoff \cite{abikoff-aug} who
proved that the mapping class group $\mathrm{Mod}(S)$ acts  properly discontinouly on $\widehat{\mathcal T}(S)$ and that  its quotient is homeomorphic
to the Deligne-Mumford compactification of the Moduli space of $S$. (Recall that the Moduli space of $S$  is the quotient of $\mathcal T(S)$ by the action of $\mathrm{Mod}(S)$.)
As such, one may view the augmented Teichm\"uller space as the ``orbifold universal cover'' of the Deligne-Mumford compactification of Moduli space.

The most concrete way to describe the augmented Teichm\"uller space is to look at the local coordinates given by extending the Fenchel-Nielsen coordinates
on Teichm\"uller space. Suppose that $P=\{\gamma_1,\ldots,\gamma_{3g-3}\}$ is a pants decomposition, i.e. a collection of disjoint simple closed curves decomposing
$S$ into $2g-2$ pairs of pants (subsurfaces homeomorphic to a twice-punctured disk). The Fenchel-Nielsen coordinates associate to each hyperbolic
surface in $\mathcal T(S)$ and each curve a positive  real co-ordinate which is the length of the associated geodesic in the surface and another real coordinate
which records the ``twist'' about the curve. (The choice of twist co-ordinate involves additional choices, but is canonical once appropriate choices are made.) This results in a real analytic
coordinate system for Teichmuller space as
$$\mathcal T(S)\cong (\mathbb R_{>0}\times \mathbb R)^{3g-3}.$$

If $C$ is any non-empty sub-collection of $P$, then there is a ``stratum at infinity'' consisting of (marked) finite area hyperbolic surfaces homeomorphic to
$S\setminus C$. One naturally, obtains Fenchel-Nielsen coordinates on this strata by looking at the curves in $P\setminus C$, so
$$\mathcal T(S\setminus C)\cong (\mathbb R_{>0}\times \mathbb R)^{\#(P\setminus C)}\cong \prod_{R\in S\setminus C}\mathcal T(R)$$
where $R$ is a component of $S\setminus C$.
If $C=\{\gamma_1\}$, then one may append the  stratum $\mathcal T(S\setminus\{\gamma_1\})$ by allowing the length coordinate associated to $\gamma_1$ to be 0 and 
forgetting the twist co-ordinate when this occurs. So,
$$\mathcal T(S)\cup\mathcal T(S\setminus\{\gamma_1\})\cong  (\mathbb R_{>0}\times \mathbb R)^{3g-4}\times  (\mathbb R_{\ge 0}\times \mathbb R)/\sim$$
where the equivalence relation is given by letting $(0,s)\sim (0,t)$ for all $s,t\in\mathbb R$.
Notice that the resulting space has the unfortunate property of failing to be locally compact.
More generally, one may append all  the strata at infinity which are Teichmuller spaces of surfaces pinched along sub-collections of $P$ at the same time, to obtain
$$\widehat{\mathcal T}^P(S)= \left((\mathbb R_{\ge 0}\times \mathbb R)/\sim\right)^{3g-3}.$$
Then, one defines  the augmented Teichm\"uller space $\widehat{\mathcal T}(S)$ to be the union of the $\widehat{\mathcal T}^P(S)$ over the collection $\mathcal P$ of
all pants decompositions of $S$,  with the obvious
identification of (marked) isometric surfaces and the topology induced by regarding the $\widehat{\mathcal T}^P(S)$ as local coordinate systems, so
$$\widehat{\mathcal T}(S)=\bigcup_{P\in\mathcal P}\widehat{\mathcal T}^P(S).$$
We regard $\mathcal T(S)=\mathcal T(S\setminus\emptyset)$ as the central stratum and each $\mathcal T(S\setminus C)$ as a stratum at infinity.
This description is  discussed more fully in Bill  Abikoff's book \cite[Sec. II.3.4]{abikoff-book}.

One may also give a more representation-theoretic viewpoint on the augmented Teichm\"uller space. One may regard an element of
$\mathcal T(S\setminus C)\cong\prod_{R\in S\setminus C}$ as a collection $\rho=\{\rho_R\}_{R\in S\setminus C}$ of (conjugacy classes of ) representations of 
$\pi_1(R)$ into $\mathsf{PSL}(2,\mathbb R)$
where $R$ is a component of $S\setminus C$ and $\rho_R$  takes each curve freely homotopic into the boundary of $R$ to a parabolic element. 
If each $\rho_i=\{(\rho_i)_{R_i}\}_{R_i\in S\setminus C_i}$ is a collection of elements of $\mathcal T(S\setminus C_i)$ where each $C_i$ is a collection of disjoint simple closed curves on $S$,
then we say that $\{\rho_i\}$ converges to $\{ (\rho_\infty)_{R_\infty}\}_{R_\infty\in S\setminus C_\infty}$ if $C_i$ is contained in $C_\infty$ for all large enough $i$,
and for each component $R_\infty$ of $S\setminus C_\infty$, the sequence $\{\rho_i|_{R_{\infty}}\}$ (which makes sense for large enough $i$) converges (up to conjugation)
to $(\rho_\infty)_{R_\infty}$. This viewpoint is explored  more fully in the larger setting of representations into $\mathsf{PSL}(2,\mathbb C)$ by Dick Canary and Pete Storm
\cite{canary-storm}.

\medskip

Howard Masur \cite{masur-wp} showed that the augmented Teichm\"uller space is homeomorphic to the metric completion of Teichm\"uller space
with the Weil-Petersson metric. Morever, the induced metric on each stratum at infinity is its Weil-Petersson metric. We will formulate his result in a manner
which is convenient for our later description of the conjectural geometric picture of the augmented Hitchin component.

\begin{theorem}{\rm (Masur \cite{masur-wp})}
Let $d_C^{WP}$ denote the Weil-Petersson metric on $\mathcal T(S\setminus C)$ where $C$ is a (possibly empty) collection of disjoint non-parallel simple closed curves.
There exists a complete metric $\hat d$ on $\widehat{\mathcal T}(S)$ such that if $\mathcal T(S\setminus C)$ is any stratum of $\widehat{\mathcal T}(S)$
then the restriction of $\hat d$ to $\mathcal T(S\setminus C)$ agrees with $d^{WP}_C$.
\end{theorem}

\medskip\noindent
{\bf Augmented Hitchin components:}
Suhyoung Choi and Bill Goldman \cite{goldman-choi} showed that if $\rho\in \mathcal H_3(S)$, then there exists a strictly convex domain $\Omega_\rho$ in
the projective plane $\mathbb{RP}^2$ so that $\rho(\pi_1(S))$ acts freely and properly discontinuously on $\Omega_\rho$, so
$X_\rho=\Omega_\rho/\rho(\pi_1(S))$ is a strictly convex projective surface. (A domain in $\mathbb{RP}^2$ is strictly convex if its closure lies
in an affine chart for $\mathbb{RP}^2$ and is strictly convex in that chart.) Bill Goldman \cite{goldman-coordinates} produced coordinates on
$\mathcal H_3(S)$ which generalize the classical Fenchel-Nielsen coordinates. 

John Loftin and Tengren Zhang \cite{loftin-zhang} use a modification of Goldman's  coordinates for $\mathcal H_3(S)$ which was developed
by Tengren Zhang \cite{zhang1}. One again begins with a pants decomposition $P=\{\gamma_1,\ldots,\gamma_{3g-3}\}$ of $S$. For each curve $\gamma_i$
there are two length coordinates, given by $\alpha_1(\rho(\gamma_i))$ and $\alpha_2(\rho(\gamma_i))$, which determine $\rho(\gamma_i)$ up to conjugacy,
and two real-valued coordinates, one of which extends the twist coordinate from the Fuchsian setting and the other of which is called the bulge coordinate.
Associated to every pair of pants in $S-P$ there are two real-valued coordinates, which are called internal coordinates. The projective structure on
each pair of pants in $S-P$ is determined by the length coordinates of its boundary curves together with its two internal coordinates. In these coordinates,
$$\mathcal H_3(S)\cong \left( (\mathbb R_{>0})^2\times\mathbb R^2\right)^{3g-3}\times \left(\mathbb R^2\right)^{2g-2}.$$

If $C$ is a non-empty subcollection of the pants decomposition $P$, then we let $\mathcal H_3(S\setminus C)$ denote the space of marked finite area projective structures on
$S\setminus C$. Then, both length coordinates associated to any curve in $C$ are $0$ and we have no twist or bulge co-ordinate associated to any curve in $C$, so
$$\mathcal H_3(S\setminus C)\cong \left( (\mathbb R_{>0})^2\times\mathbb R^2\right)^{\#(P\setminus C)}\times \left(\mathbb R^2\right)^{2g-2}\cong\prod_{R\in S\setminus C}\mathcal H_3(R)$$
where $\mathcal H_3(R)=\mathcal H_3(\Gamma_R)$ where $\mathbb H^2/\Gamma_R$ is a finite area hyperbolic surface homeomorphic to $R$.
(Ludovic Marquis \cite{marquis-deformation} also describes a parametrization of $\mathcal H_3(S\setminus C)$.)
If $C=\{\gamma_1\}$,
we append $\mathcal H_3(S\setminus \{\gamma_1\})$ to $\mathcal H_3(S)$ as a ``stratum at infinity'' much as we did in the Teichm\"uller setting 
by allowing the additional length coordinate pair $(0,0)$ for the curve $\gamma_1$  but ignoring the twist and bulge coordinates whenever this occurs, so
$$\mathcal H_3(S)\cup\mathcal H_3(S\setminus\{\gamma_1\})\cong  \left( (\mathbb R_{>0})^2\times\mathbb R^2\right)^{3g-4}\times  
\Big( (\mathbb R_{\ge0})^2\times\mathbb R^2\Big)/\sim\times \left(\mathbb R^2\right)^{2g-2}$$
where $(0,0,a,b)\sim(0,0,c,d)$ for all $a,b,c,d\in\mathbb R$. Again, the resulting space fails to be locally compact.
If we let $\widehat{\mathcal H}_3^P(S)$ be the result of appending $\mathcal H_3(S\setminus C)$ for all non-empty subcollections $C$ of $P$, then the
resulting space is parameterized as
$$\widehat{\mathcal H}^P_3(S)\cong 
\Big( \big( (\mathbb R_{\ge0})^2\times\mathbb R^2\big)/\sim\Big)^{3g-3}\times \left(\mathbb R^2\right)^{2g-2}.$$
Then, one defines the augmented Hitchin component $\widehat{\mathcal H}_3(S)$ to be the union of the $\widehat{\mathcal H}_3^P(S)$ over the collection $\mathcal P$ of
all pants decompositions of $S$,  with the obvious
identification of (marked)  convex projective surfaces and the topology induced by regarding the $\widehat{\mathcal H}_3^P(S)$ as local coordinate systems, so
$$\widehat{\mathcal H}_3(S)=\bigcup_{P\in\mathcal P}\widehat{\mathcal H}_3^P(S).$$
Notice that if $\mathcal C$ is the set  of (isotopy classes of) non-empty collections of disjoint non-parallel curves on $S$, then
$$\widehat{\mathcal H}_3(S)=\mathcal H_3(S)\bigsqcup_{C\in\mathcal C} \mathcal H_3(S\setminus C)$$
and we regard $\mathcal H_3(S)$ as the central stratum and each $\mathcal H_3(S\setminus C)$ as a stratum at infinity.

Loftin and Zhang consider a larger augmented Hitchin component where the ``pinched'' curves are not required to map to unipotent elements. 
This space is natural from various viewpoints. However, we expect, with very little evidence, that our smaller augmented Hitchin component arises
as the metric completion of the Hitchin component with a pressure metric. John Loftin \cite{loftin} earlier gave a more  analytic description of the augmented Hitchin
component from the point of view of the parametrization of $\mathcal H_3(S)$ as the bundle of cubic differentials over Teichm\"uller space. This more analytic viewpoint is also 
useful in the study of the augmented Hitchin component, but we will not discuss it further here.

\medskip\noindent
{\bf The Geometric Picture:}
We are now ready to give a description of our conjectural geometric picture. If $\phi\in\Delta$, then there is a pressure form $\mathbb P^\phi$
on each ``stratum at infinity'' of $\widehat{\mathcal H}_3(S)$ and on the central stratum $\mathcal H_3(S)$. This gives rise to a path pseudo-metric $d^\phi$ on each stratum.
If $\phi$ is $\omega_1$ or $\alpha_1$, then $d_\phi$ is an analytic Riemannian metric on each stratum, while if $\phi=\omega_H$ then it is a
path metric on each stratum. One way to state our conjecture is the following.

\medskip\noindent
{\bf Conjecture:} {\em If $\phi$ is $\omega_1$, $\alpha_1$ or $\omega_H$, then $d^\phi$ extends to a complete metric $\hat d^\phi$ on $\widehat{\mathcal H}_3(S)$.
In particular, $\widehat{\mathcal H}_3(S)$ is the metric completion of $\mathcal H_3(S)$ with the pressure metric $d^\phi$.}

\medskip

Even though  the metric $d^{\omega_H}$ is only a path metric on each stratum, it may be the most natural metric to investigate. If $\phi\in\mathcal H_3(S)$,
then the Hilbert metric on $\Omega_\phi$ is Gromov hyperbolic and induces a Finsler metric on $X_\rho=\Omega_\rho/\rho(\pi_1(S))$ and $\ell^{\omega_H}(\rho(\gamma))$ is
the length of the unique closed geodesic on $X_\rho$ in the free homotopy class of $\gamma$, see Marquis \cite{marquis-finite}.

\medskip

One expects that the conjectural picture described above should generalize to all values of $d$. In this case, one should be able to replace the Zhang-modified
Goldman coordinates, with the coordinates described by Fock-Goncharov \cite{fock-goncharov} or Bonahon-Dreyer \cite{bonahon-dreyer}. The representation
theoretic viewpoint also generalizes, where one replaces the assumption that the images of boundary components of complementary surfaces map to parabolic
elements with the assumption that the restriction to each complementary surface $R$ is a Hitchin representation of a finite area uniformization of $R$.

\end{document}